\documentclass{article}
\usepackage[utf8]{inputenc}
\usepackage{amsmath}
\usepackage{amsfonts}
\usepackage{amssymb}
\usepackage{amsthm}
\usepackage{nicefrac}

\usepackage{microtype}
\usepackage[dvipsnames]{xcolor} %
\usepackage{tikz}
\usepackage{pgfplots}
\usetikzlibrary{calc,decorations.markings,positioning,arrows}
\pgfplotsset{compat=1.17}

\theoremstyle{plain}

\newtheorem{theorem}{Theorem}

\newtheorem{remark}{Remark}

\theoremstyle{definition}

\usepackage[export]{adjustbox}
\usepackage{subfig}
\usepackage{graphicx}

\usepackage{booktabs}
\usepackage{multirow}
\usepackage{todonotes}

\usepackage[]{algorithm2e}
\usepackage{hyperref}
\usepackage{matlab-prettifier}
\usepackage{authblk}

\makeatletter
\newcommand{\subjclass}[2][2020]{%
	\let\@oldtitle\@title%
	\gdef\@title{\@oldtitle\footnotetext{#1 \emph{Mathematics subject classification.} #2}}%
}
\newcommand{\keywords}[1]{%
	\let\@@oldtitle\@title%
	\gdef\@title{\@@oldtitle\footnotetext{\emph{Key words and phrases.} #1.}}%
}
\makeatother

\author{Lidia Aceto}\affil{Dipartimento di Scienze e Innovazione Tecnologica, Universit\`{a} del Piemonte Orientale,\\viale T. Michel 11, 15121 Alessandria, Italy,\\E-mail: lidia.aceto@uniupo.it}
\author{Fabio Durastante}\affil{Dipartimento di Matematica, Universit\`{a} di Pisa,\\via F. Buonarroti 1/C, 56127 Pisa, Italy,\\E-mail: fabio.durastante@unipi.it}
\title{Efficient computation of the Wright function and its applications to fractional diffusion-wave equations}
\date{\today}

\begin{document}

	\maketitle

    \begin{abstract}
	In this article, we deal with the efficient computation of the Wright function in the cases of interest for the expression of solutions of some fractional differential equations. The proposed algorithm is based on the inversion of the Laplace transform of a particular expression of the Wright function for which we discuss in detail the error analysis. We also present a code package that implements the algorithm proposed here in different programming languages. The analysis and implementation are accompanied by an extensive set of numerical experiments that validate both the theoretical estimates of the error and the applicability of the proposed method for representing the solutions of fractional differential equations.
\end{abstract}

\section{Introduction}
The Wright function is a generalization of the exponential function defined by the convergent series in the whole complex plane
\begin{equation} \label{wrightfun}
W_{\lambda,\mu}(z) := \sum_{n=0}^\infty \frac{z^n}{n! \, \Gamma(\lambda n + \mu)}, \quad \lambda >-1, \, \mu \in \mathbb{C},
\end{equation}
where $\Gamma$ is the Euler gamma function. It was initially introduced by E.M. Wright for $\lambda \ge 0$ in the framework of the asymptotic theory of partitions \cite{Wright}. The same author then studied the case $\lambda \in (-1,0),$ which is now referred to in the literature as  {\em Wright function of the second kind} (WF2K)  \cite{Wright1940}. 

Although several representations of the Wright function have been introduced such as integral representation and asymptotic expansions and many of its analytical properties have already been well-studied (see, e.g. \cite{Stankovic,Luchko+2019,MainardiConsiglio,mainardi2010m}), its numerical evaluation is still an active research area. In this paper, we discuss the numerical evaluation of WF2K, since this is the most interesting case for applications. Indeed, this function plays an important role in describing non-Gaussian deterministic and stochastic processes and the transition from sub-diffusion processes to wave propagation. Noteworthy cases of WF2Ks are the functions 
\begin{equation}\label{eq:wrightasspecialfunction}
F_\nu(z) :=  W_{-\nu,0}(-z), \qquad \quad M_\nu(z) :=  W_{-\nu,1-\nu}(-z), \quad \nu \in (0,1),
\end{equation}
related via the formula  $F_\nu(z)= \nu z M_\nu(z).$ 
They were introduced by Mainardi in the 90s and are called {\em auxiliary functions} by virtue of their roles in solving the {signalling} problem and the Cauchy value problem, respectively, for the time-fractional diffusion-wave equation (see \cite{mainardi2020tutorial} and references therein). Except for some special cases where WF2K can be represented in terms of other elementary and special functions, such as 
\begin{equation}\label{eq:special-cases}
M_0(z)= \exp(-z), \qquad \qquad M_{1/2}(z)=1/\sqrt{\pi} \exp(-z^2/4), \qquad \qquad M_{1/3}(z)=3^{2/3} \mbox{Ai}(z/3^{1/3}),
\end{equation}
for $\mbox{Ai}$ the Airy function~\cite[Chapter~9]{NIST:DLMF}, most programming languages do not provide built-in functions for WF2K. In this paper we try to fill this gap by considering the numerical evaluation of the following function
\begin{equation} \label{eq:flm}
f_{\lambda,\mu}(t;x):= t^{\mu-1} W_{\lambda, \mu} (- |x| t^{\lambda}), \qquad  t\ge 0,\, x \in \mathbb{R}, \quad \lambda \in (-1,0),\, \mu \in \mathbb{C},
\end{equation}
obtained from \eqref{wrightfun} by setting $z=- |x| t^{\lambda}.$ 
As an example, in Figure~\ref{fig:function_examples} we show the plots of some of these functions.
\begin{figure}[htbp]
	\centering
	\input{function_examples.tikz}
	\caption{Plots of the function $f_{-\nu,1-\nu}(1;x)$ for $\nu=\nicefrac{1}{2}$, $\nicefrac{3}{8}$, $\nicefrac{1}{4}$, and $\nicefrac{1}{8}$.}
	\label{fig:function_examples}
\end{figure}

Starting from the integral representation of the Wright function several algorithms for its numerical computation have already been proposed in \cite{LuchkoFCAA2008,luchko2010wright}. The authors state that these could serve as the basis for the creation of a programming package for the numerical evaluation of the Wright function. However, the actual implementation is not straightforward as it seems. The proposed integrals have oscillatory behavior and selecting efficient quadrature formulas or a reliable truncation seems to be challenging.   

In this paper, we consider a different approach based on the inversion of the Laplace transform in which a trapezoidal rule is applied on a parabolic contour. This approach solves the difficulties encountered in adopting the strategies discussed in~\cite{LuchkoFCAA2008,luchko2010wright}. 
It is worth mentioning that methods of this type have already been successfully applied in \cite{DingfelderWeideman,GarrappaPopolizioMLRealLine} for the computation of Mittag-Leffler functions on the real line, { in using a solve-then-discretize approach for certain fractional differential equations~\cite{MR4371136} and, more generally, for the computation of actions of semigroups with certified error bounds~\cite{MR4379629}}.

The paper is organized as follows. 
In section \ref{sec:inversion_of_the_Laplace_transform} we first recall some generalities on the computation of Laplace-type integrals. Then we focus on our case of interest and delineate the relative error analysis leveraging corresponding results for the Mittag-Leffler function from~\cite{GarrappaPopolizioMLRealLine}. In section~\ref{sec:time-fractional-evolution-equation} we look at several fractional differential problems whose solution can be expressed in terms of the Wright function to substantiate some possible use cases of the algorithm presented here. In section \ref{sec:numerexp} we present some numerical experiments which confirm our analysis from section~\ref{sec:inversion_of_the_Laplace_transform}. Finally, section \ref{sec:conclusion} summarizes the main conclusions of this work.
\section{Inversion of the Laplace transform}\label{sec:inversion_of_the_Laplace_transform}

\begin{figure}[htbp]
	\centering
	\begin{tikzpicture}[scale=0.7]
	\begin{axis}[
	axis x line=center,
	axis y line=center,
	xlabel={$\mathbb{R}$},
	ylabel={$i\mathbb{R}$},
	ticks=none,
	xlabel style={below right},
	ylabel style={above left},
	xmin=-4.5,
	xmax=4.5,
	ymin=-3.5,
	ymax=3.5]
	\draw[red, dashed, line width=1.25pt, domain={-2.5}:{3}, samples=500] plot[smooth]({-\x*\x +1},\x );
	\filldraw [black] (0.00, 0.00) circle(3pt);
	\draw[-,line width=1.7pt] (0,0) to (-4.5,0);
	\fill [color=red] (1,0) circle (2pt);
	\draw[color=black] (1.6,0.35) node {$(\gamma,0)$};
	\end{axis}\hfil
	\end{tikzpicture}
	\begin{tikzpicture}[scale=0.7]
	\pgfmathsetmacro{\N}{10}
	\pgfmathsetmacro{\h}{3/\N}
	\begin{axis}[
	axis x line=center,
	axis y line=center,
	xlabel={$\mathbb{R}$},
	ylabel={$i\mathbb{R}$},
	ticks=none,
	xlabel style={below right},
	ylabel style={above left},
	xmin=-4.5,
	xmax=4.5,
	ymin=-3.5,
	ymax=3.5]
	\draw[red, dashed, line width=0.25pt, domain={-2.5}:{3}, samples=500] plot[smooth]({-\x*\x +1},\x );
	\addplot[only marks,mark=x,mark size=4pt] coordinates {
		(	   -3.0000	,	 -2.0000	)
		(	   -1.4198	,	 -1.5556	)
		(	   -0.2346	,	 -1.1111	)
		(	    0.5556	,	 -0.6667	)
		(	    0.9506	,	 -0.2222	)
		(	    0.9506	,	  0.2222	)
		(	    0.5556	,	  0.6667	)
		(	   -0.2346	,	  1.1111	)
		(	   -1.4198	,	  1.5556	)
		(	   -3.0000	,	  2.0000	)
	};
	\end{axis}
	\end{tikzpicture}
	\caption{Parabolic contour {\protect\raisebox{0.25em}{\protect\tikz{\protect\draw[red, -, dashed, line width=1.2pt]  (0,0) to (0.5,0);}}}, singularity {\protect\raisebox{0.15em}{\protect\tikz{\protect\filldraw [black] (0.00, 0.00) circle(2pt);}}}, quadrature nodes $\times$, and branch cut {\protect\raisebox{0.25em}{\protect\tikz{\protect\draw[black, -, line width=1.2pt]  (0,0) to (0.5,0);}}} of the $F_{\lambda,\mu}(s;x)$  function given in \eqref{eq:Flm}.}
	\label{fig:parabolic_contour}
\end{figure}
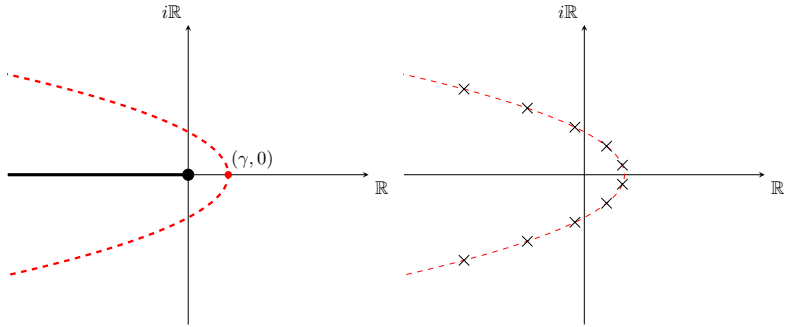

The aim of this section is to show that an effective way to compute $f_{\lambda,\mu}(t;x)$ given in \eqref{eq:flm} is based on the numerical inversion of the Laplace transform
\begin{equation} \label{eq:inverLapl}
f_{\lambda,\mu}(t;x) = \frac{1}{2\pi i} \int_{\mathcal{C}} e^{st} F_{\lambda,\mu}(s;x)\,{\rm d}s,
\end{equation}
where (see, e.g. \cite[Eq. (11)]{Stankovic})
\begin{equation} \label{eq:Flm}
F_{\lambda,\mu}(s;x) = s^{-\mu} e^{-|x| s^{-\lambda}}.
\end{equation}
Here the contour $\mathcal{C}$ is a suitable deformation of the Bromwich line that must be fixed in a region of the complex plane in which $F_{\lambda,\mu}(s;x)$ is analytic. 
Considering that the only singularity of $F_{\lambda,\mu}(s;x)$ lies in the branch point $s=0,$ we can choose the simplest contour so far discussed in the literature, namely the parabola~\cite{WeidemanTrefethenParHyp}. This contour is described by the equation 
\[
z(u)= \gamma (i u +1)^2, \qquad -\infty < u < +\infty,
\]
where $\gamma$ is a parameter to be taken positive so that the parabola encloses the entire negative real axis; see Figure \ref{fig:parabolic_contour}.
Using this parabolic-shaped contour $z(u)$ we can write the integral \eqref{eq:inverLapl} as
\begin{eqnarray} \label{eq:intzu}
f_{\lambda,\mu}(t;x) &=& \frac{1}{2\pi i} \int_{-\infty}^{+\infty}  e^{z(u) t} F_{\lambda,\mu}(z(u);x) z'(u)\,{\rm d}u \nonumber\\
&=& \frac{1}{2\pi i} \int_{-\infty}^{+\infty}   g_{\lambda,\mu}(u)\,{\rm d}u,
\end{eqnarray}
where
\begin{equation} \label{eq:glm}
g_{\lambda,\mu}(u):=e^{z(u) t}  (z(u))^{-\mu} e^{-|x| (z(u))^{-\lambda}} z'(u).
\end{equation}
By approximating this integral by means of the finite trapezoidal rule with step-size $h$ we obtain
\begin{equation} \label{eq:trapN}
f^{[h,N]}_{\lambda,\mu}(t;x) = \frac{h}{2\pi i} \sum_{k=-N}^N  g_{\lambda,\mu}(u_k), \qquad u_k=kh.
\end{equation}
We observe that for the symmetry of the contour with respect to the real axis, only the quadrature nodes in the upper (lower) half-plane can be  considered; this leads to a halving of the computational cost. In addition, for the given $t$ and $x$, the three parameters $\gamma, h$ and $N$ 
occurring in  (\ref{eq:trapN}) must be set to minimize the error 
\begin{equation} \label{eq:errabs}
E^{[h,N]}(t;x):=\left| f_{\lambda,\mu}(t;x) - f^{[h,N]}_{\lambda,\mu}(t;x) \right|.
\end{equation}
Before carrying out a detailed error analysis, for the reader's convenience, we recall some basic facts about the well-known error estimates for the trapezoidal rule that will be useful for such analysis.
\subsection{Error estimates for the trapezoidal rule} \label{sec:trap}
Consider the absolutely convergent integral
\begin{equation*}
I=  \int_{-\infty}^{+\infty} g(u)\,{\rm d}u,
\end{equation*}
and its infinite and finite trapezoidal approximations
\begin{equation*}
I^{[h]}= h \sum_{k=-\infty}^{+\infty}  g(kh), \qquad  I^{[h,N]}=h \sum_{k=-N}^N  g(kh).
\end{equation*}
For the error we have
\begin{equation} \label{eq:dueerrori}
\left|I-I^{[h,N]} \right|  \le \mathcal{E}_{D} +\mathcal{E}_{T},
\end{equation}
where 
\[
\mathcal{E}_{D}  =  \left|I-I^{[h]} \right|, \qquad 
\mathcal{E}_{T} = \left|I^{[h]} - I^{[h,N]} \right|.
\]
The quantities $\mathcal{E}_{D}$ and $\mathcal{E}_{T}$ are often referred to as the discretization error and the truncation error, respectively. 

As for the discretization error, when $g(u)$ is a complex-valued function we can use the following theorem.
\begin{theorem}\cite[Theorem 2.1]{WeidemanTrefethenParHyp} \label{teo:DE}
	Let $w=u+iv,$ with $u$ and $v$ real. Suppose g(w) is analytic in the strip $-d<v<c,$ for some $c>0,d>0,$ with $g(w)\rightarrow 0$ uniformly as $|w| \rightarrow +\infty$ in that strip. Suppose further that for some $M_+>0, M_->0$ the function $g(w)$ satisfies
	\begin{equation*}
	\int_{-\infty}^{+ \infty} \left| g(u+ir) \right|\,{\rm d}u \le M_+, \qquad  \int_{-\infty}^{+ \infty} \left| g(u-is) \right|\,{\rm d}u \le M_-,
	\end{equation*}
	for all $0<r<c, 0<s<d.$ Then,
	\[
	\mathcal{E}_{D} \le   \mathcal{E}_{D_+}  +   \mathcal{E}_{D_-}, 
	\]
	where 
	\[
	\mathcal{E}_{D_+} = \frac{M_+}{e^{2\pi c/h} -1}, \qquad  \mathcal{E}_{D_-}= \frac{M_-}{e^{2\pi d/h} -1}.
	\]
\end{theorem}

\begin{remark}
	When  the contribution of $M_+$ and $M_-$ is negligible then the estimates $\mathcal{E}_{D_+} \approx e^{-2\pi c/h}$ and $ \mathcal{E}_{D_-}\approx e^{-2\pi d/h}$ are sufficient to obtain a satisfactory error analysis. Otherwise, these terms would have to be balanced with the decay of the exponential to get good estimates. For example, this is the situation  encountered when dealing with the computation of Mittag-Leffler functions on the real line~\cite{GarrappaPopolizioMLRealLine}{, and when dealing with the numerical approximation of strongly continuous semigroups on infinite-dimensional Hilbert spaces~\cite{MR4379629}.}
\end{remark}
As for the truncation error, if $g(u)$ decays rapidly as $u \rightarrow \pm \infty,$ we have 
\[
\mathcal{E}_{T} = \mathcal{O}(|g(hN)|), \quad N \rightarrow +\infty,
\]
that is,  the truncation error can be approximated by the last term retained of the trapezoidal sum.

Following %
\cite{weideman2010improved}, to make a more complete error analysis, we also consider the roundoff errors. At this aim we introduce the term
\begin{equation*} 
I^{[h,N,\epsilon]}= h \sum_{k=-N}^{N}  g(kh) (1+\epsilon_k),
\end{equation*}
where $\epsilon_k$  are the relative errors in the computed function values $g(kh)$ and they all satisfy $|\epsilon_k|\le \epsilon,$ the precision machine. The formula (\ref{eq:dueerrori}) is then extended to
\begin{equation*} %
\left|I-I^{[h,N]} \right|  \le \mathcal{E}_{D} +\mathcal{E}_{T} + \mathcal{E}_{R},
\end{equation*}
where 
\[
\mathcal{E}_{R}  = \left|I^{[h,N]} - I^{[h,N,\epsilon]} \right|.
\]
Obviously,
\begin{equation} \label{eq:RE}
\mathcal{E}_{R}  \le \epsilon  h \sum_{k=-N}^{N}    \left| g(kh)\right|.
\end{equation}
\subsection{Selecting the quadrature nodes}\label{sec:selecting_the_quadrature_nodes}
We use here the results from section~\ref{sec:trap} to perform the error analysis for~\eqref{eq:trapN}. Denoting by (see \eqref{eq:glm})
\begin{equation*} %
f^{[h]}_{\lambda,\mu}(t;x) = \frac{h}{2\pi i} \sum_{k=-\infty}^{+\infty} g_{\lambda,\mu}(u_k), \qquad u_k=kh,
\end{equation*}
the infinite trapezoidal approximation to \eqref{eq:intzu}, we can bound the error in \eqref{eq:errabs} by means of the discretization error and the truncation error which in this case are defined by
\[
E^{[h]}_D (t;x) = \left|f_{\lambda,\mu}(t;x) - f^{[h]}_{\lambda,\mu}(t;x)\right|, \quad E^{[h,N]}_T (t;x) = \left|f^{[h]}_{\lambda,\mu}(t;x)-f^{[h,N]}_{\lambda,\mu}(t;x)\right|,
\]
namely,
\begin{equation} \label{eq:pez1}
E^{[h,N]}(t;x) \le E^{[h]}_D (t;x)+ E^{[h,N]}_T (t;x).    
\end{equation}
Furthermore, based on the results of Theorem \ref{teo:DE} we can use a non-symmetric strip of analyticity and therefore we also have that
\begin{equation} \label{eq:pez2}
E^{[h]}_D (t;x) \le E^{[h]}_{D_+} (t;x) + E^{[h]}_{D_-} (t;x),
\end{equation}
where the quantities $E^{[h]}_{D_+} (t;x) $ and $E^{[h]}_{D_-} (t;x) $ depend on the behavior of $g_{\lambda,\mu}$ on the upper and  on the lower half-plane, respectively. In order to determine an estimate of these two quantities, we begin by studying the behavior of $|g_{\lambda,\mu}(u+ir)|$ for  $0<r < c<1.$ Following the results reported in the proof of Theorem 2 given in \cite{GarrappaPopolizioMLRealLine} it is immediate to verify that 
\begin{align*}  
\left|e^{z(u+ir) t}\right|  & = e^{\gamma (1-r)^2 t} e^{-\gamma u^2 t} \le e^{\gamma (1-u^2) t}\\
\left| e^{-|x| (z(u+ir))^{-\lambda}}\right| & = e^{-|x| \left|z(u+ir)^{-\lambda}\right| } = e^{-|x|  \gamma^{-\lambda} ((1-r)^2 -u^2)^{-\lambda}}  \le 1 \\
\left|z'(u+ir)\right| &= 2 \gamma \sqrt{(1-r)^2 +u^2}.
\end{align*}
Concerning the term $\left|(z(u+ir))^{-\mu}\right|$ we observe that $(z(u+ir))^{-\mu} =e^{\log_{\mathbb{C}} (z(u+ir))^{-\mu}}$ and then
\begin{align*}
\left|(z(u+ir))^{-\mu} \right| & = e^{ \Re(-\mu \log_{\mathbb{C}} (z(u+ir)))} \\
& = e^{ -\Re(\mu) \log \left|z(u+ir) \right| +\mathfrak{Im}(\mu) \mathrm{arg}(z(u+ir))} \\
& = \left( \gamma \sqrt{(1-r)^2+u^2}\right)^{-\Re(\mu) }  e^{\mathfrak{Im}(\mu) \mathrm{arg}(z(u+ir))},
\end{align*}
where $\log(\cdot)$ is the natural logarithm, and $\log_{\mathbb{C}}$ is the complex logarithm. Consequently,  from~\eqref{eq:glm} we have
\begin{align*}  
\int_{-\infty}^{+ \infty}   \left|g_{\lambda,\mu}(u+ir)\right|  \,{\rm d}u &  \le  M_1  \int_{-\infty}^{+ \infty}  e^{ - \gamma u^2  t}  ((1-r)^2 +u^2)^{ \frac{1-\Re(\mu)}{2}}   \,{\rm d}u,
\end{align*}
where {$M_1= 2 \gamma^{1-\Re (\mu)} e^{\gamma t +2\pi\mathfrak{Im}(\mu)}$}. Using \cite[Lemma 1]{GarrappaPopolizioMLRealLine}  we obtain 
\begin{align*}
\int_{-\infty}^{+ \infty}   \left|g_{\lambda,\mu}(u+ir)\right|  \,{\rm d}u &  \le  M  (1-r)^{2-\Re (\mu)} \Psi\left(\frac{1}{2}, 2-\frac{\Re (\mu)}{2},  \gamma t (1-r)^2\right), \quad \text{ for } \gamma t (1-r)^2 \rightarrow 0, 
\end{align*}
where $M= \sqrt{\pi}M_1$ is a positive constant which does not depend on $r$ and $\Psi$ is the confluent hypergeometric function {of}  the second kind; see, e.g., \cite[\S 13.2]{NIST:DLMF}. 
Therefore, from Theorem 3 given in \cite{GarrappaPopolizioMLRealLine} we can deduce the following result.
\begin{theorem}
	Let $0<r < c <1.$ Then,
	\begin{align*}
	\int_{-\infty}^{+ \infty}   \left|g_{\lambda,\mu}(u+ir)\right|  du & \le  {\hat{M}} \, w(c),
	\end{align*}
	where ${\hat{M}}$ is a suitable  positive constant that depends on $\mu$ but not on $c$ and
	\begin{align} \label{eq:wc}
	w(c) = \left\{\begin{array}{lll}
	1 & \text{ if }   \Re(\mu) < 2 \\
	-\log(\gamma t (1-c)^2)  & \text{ if }     \Re(\mu) =2, & 1 - \frac{1}{\sqrt{\gamma t}} < c < 1, \\
	(1-c)^{2-\Re (\mu)} &  \text{ if } \Re(\mu) >2,  
	\end{array}
	\right.
	\end{align}
	for $\log(\cdot)$ the natural logarithm.
\end{theorem}
Using this result in Theorem~\ref{teo:DE} we obtain
\begin{align} \label{eq:dp}
E^{[h]}_{D_+} (t;x) &= \mathcal{O}(w(c) e^{-2\pi c/h}).
\end{align}
Hence, from~\eqref{eq:pez1}-\eqref{eq:pez2}, an error $E^{[h,N]}(t;x)$ proportional to a given tolerance $\varepsilon >0$ is obtained after balancing $E^{[h]}_{D_+} (t;x)$ with the standard estimates 
\begin{equation}\label{eq:truncation_and_discretization_errors}
E^{[h]}_{D_-} (t;x) = \mathcal{O}(e^{-\pi^2/(\gamma t h^2)+2\pi /h}), \qquad E^{[h,N]}_T (t;x) = \mathcal{O}\left( e^{\gamma t (1- (hN)^2)}\right), \qquad h \rightarrow 0,
\end{equation}
provided in \cite{WeidemanTrefethenParHyp}. 
However, as already mentioned at the end of section \ref{sec:trap}, a more accurate analysis also takes into account the roundoff error. The following estimate applies to this error {(substitute~\eqref{eq:glm} in \eqref{eq:RE} and use the symmetry to half the domain of the integral)%
	\begin{align}
	E^{[h,N]}_{R} (t;x) & \le \frac{\epsilon h}{2\pi} e^{\gamma t} \sum_{k=-N}^{N} \left| (z(u_k))^{-\mu} e^{-|x| (z(u_k))^{-\lambda}} z'(u_k)\right| \nonumber\\ & \le  \frac{\epsilon h}{2\pi} e^{\gamma t} \sum_{k=-N}^{N} \left| (z(u_k))^{-\mu} z'(u_k)\right| \nonumber \\
	&\approx \epsilon \frac{2 \gamma}{\pi} e^{\gamma t} \int_{0}^{Nh} \left(u^2+1\right)^{\frac{1-\mathfrak{R}(\mu)}{2}}\, ds.\label{eq:hypergeom}
	\end{align}
	Assuming that $h = \mathcal{O}(N^{-2}),$ we can compute the integral in~\eqref{eq:hypergeom} using the integral representation of the Gauss hypergeometric function ${}_2F_1(a,b;c;z)$ from \cite[Eq.s (15.2.1), (15.2.2), and (15.6.1)]{NIST:DLMF} with $a = \nicefrac{(-1+\Re(\mu))}{2}$, $b = \nicefrac{1}{2}$, $c=\nicefrac{3}{2}$, $z = -N^{-2}$ finding
	\[
	\int_0^{\frac{1}{N}} \left(u^2+1\right)^{\frac{1-\mathfrak{R}(\mu)}{2}} \, du = \frac{\, _2F_1\left(\frac{1}{2},\frac{-1+\mathfrak{R}(\mu)}{2};\frac{3}{2};-\frac{1}{N^2}\right)}{N} \overset{N \rightarrow +\infty}{\longrightarrow}0.
	\]
	Therefore, the integral in relation~\eqref{eq:hypergeom} can be neglected together with the scaling term  $2\gamma/\pi$ giving us}
\begin{equation}\label{eq:ERHN}
E^{[h,N]}_{R} (t;x) \approx \epsilon e^{\gamma t}.
\end{equation}

Thus,  the expression for the optimal parameters $\gamma, h$ and $N$ 
is determined by requesting that
\begin{equation} \label{eq:balance}
E^{[h]}_{D_+} (t;x) \approx E^{[h]}_{D_-} (t;x) \approx 
E^{[h,N]}_T (t;x) \approx E^{[h,N]}_{R} (t;x) 
\end{equation}
asymptotically as $h \rightarrow 0.$ First of all, in order to simplify the notation,  we denote by (see (\ref{eq:wc}))
\begin{equation} \label{eq:valcsi}
\xi = 2-\frac{h}{\pi c} \log(w(c)).
\end{equation}
In this way the balancing of the errors \eqref{eq:balance} leads to the following equations (see \eqref{eq:dp}, \eqref{eq:truncation_and_discretization_errors}, and \eqref{eq:ERHN})
\begin{equation*}%
-\frac{\xi \pi c}{h} = -\frac{\pi^2}{\gamma t h^2} +\frac{2\pi}{h}=\gamma t (1-(hN)^2) = \gamma t -\ell,
\end{equation*}
where $\ell =-\log \epsilon.$ Solving them, we get the expressions for the optimal parameters:
\begin{equation} \label{eq:oppar}
N=\frac{(2+\xi c)\ell}{(1+\xi c)\pi},  \qquad h=\frac{(2+\xi c) \ell}{\pi N^2}, \qquad   \gamma=\frac{\pi^2 N^2}{(2+\xi c)^2 t \ell}.
\end{equation}
Therefore, an error $E^{[h,N]}(t;x)\approx \varepsilon$  is obtained by requiring $e^{\log(\varepsilon)} \approx e^{- \nicefrac{\xi \pi c }{h}}$. Then, using the value of $h$ given in~\eqref{eq:oppar}, we select a number of nodes 
\begin{equation}\label{eq:valueofN}
N = \frac{\sqrt{\ell (-\log \varepsilon)}}{\pi}\sqrt{1 +\frac{2}{\xi c}}, 
\end{equation}
or, more precisely, the greatest integer less than or equal to this quantity. 

Let us observe now that $N$ is a function of $c$ which depends also on $\mu$ as we may deduce from the expression of $\xi$ in~\eqref{eq:valcsi} and by taking into account~\eqref{eq:wc}. We distinguish now between the three different cases to uncover the determination for the discretization parameters $N$, $h$ and $\gamma$. When $\Re(\mu) < 2$ we have $\xi=2$. Furthermore, as regards the value of $c$ we can determine it with the aim of minimizing the number of nodes or, equivalently, the computational cost. Since $N$ is a decreasing function in $c$
and $\sqrt{1+1/c}\rightarrow \sqrt{2} $  as $c\rightarrow 1$ 
we choose
\begin{equation} \label{eq:parmin2}
N=\left \lfloor{ \frac{\sqrt{2 \ell (-\log \varepsilon)}}{\pi}}\right \rfloor,  \qquad h=\frac{4 \ell}{\pi N^2}, \qquad \gamma=\frac{\pi^2 N^2}{16 t \ell},
\end{equation}
where $\lfloor{\cdot}\rfloor$ denotes the floor function. 

When $\Re(\mu) \ge 2,$ the value of  $c$ must be suitably fixed considering the effects of the term $w(c)$ (see~\eqref{eq:wc}). 
We start by studying in more detail what happens when $\Re(\mu)>2.$
Using the relation
\[
-\frac{\xi \pi c}{h} \approx {\log \varepsilon}
\]
in the expression of $\xi$ for this case
\[
\xi = 2-\frac{h}{\pi c} (2-\Re(\mu)) \log(1-c),
\]
we get
\begin{equation}\label{csim2}
\frac{2}{\xi} =  1+  \frac{2-\Re(\mu)}{(-\log \varepsilon)}  \log(1-c).
\end{equation}
Consequently, substituting this quantity in~\eqref{eq:valueofN} we have
\[
N(c)=   \frac{\sqrt{\ell (-\log \varepsilon)}}{\pi}\sqrt{ 1 +
	\frac{1}{c}\left( 1+  \frac{2-\Re(\mu)}{(-\log \varepsilon)}  \log(1-c) \right)}, \qquad c \in (0,1),
\]
which is a positive function that admits at least a minimum since $N(c) \rightarrow +\infty$ for both  $c\rightarrow 0$ and $c\rightarrow 1.$ We refer to it as $\bar c$. Then, $N=\lfloor{N({\bar c})}\rfloor.$ Using  $\bar c$ in (\ref{csim2}) we can compute the corresponding value of $\xi.$ Finally, the so obtained values of $c, \xi$ and $N$ can be used in (\ref{eq:oppar}) to compute $h$ and $\gamma.$

When $\Re(\mu) = 2$ using arguments similar to those just considered in the case $\Re(\mu) > 2$ and taking into account that
\[
-\frac{\xi \pi c}{h} = \gamma t - \ell \approx {\log \varepsilon},
\]
we obtain
\[
\frac{2}{\xi} =  1+  \frac{1}{(-\log \varepsilon)}  \log\left(-\log \left( (\ell+\log \varepsilon) (1-c)^2 \right)\right).
\]
This implies that 
\[
N(c)=   \frac{\sqrt{\ell (-\log \varepsilon)}}{\pi}\sqrt{ 1 +
	\frac{1}{c}\left( 1+  \frac{1}{(-\log \varepsilon)}  \log\left( -\log \left( (\ell+\log \varepsilon) (1-c)^2 \right)\right) \right)}, \qquad c \in (1 - \nicefrac{1}{\sqrt{\ell+\log(\varepsilon)}},1).
\]
Proceeding as in the case $\Re(\mu) > 2$ we can determine the value of the free parameters $N,h,$ and $\gamma.$

\begin{remark}
	To determine the parameters $h$, $\gamma$ in~\eqref{eq:oppar}, we evaluate the minimum of the function $N(c)$, for the different ranges of $\mu$, using  Brent's method. This is implemented in \textsc{MATLAB}\textsuperscript{\textregistered} by the routine \lstinline[style=Matlab-editor]{fminbnd} for constrained minimization, see Figure~\ref{fig:optimalxi}. 
		\begin{figure}[htbp]
	\begin{tabular}{p{0.48\columnwidth}p{0.46\columnwidth}}
		\begin{minipage}{.48\columnwidth}
\begin{lstlisting}[
style=Matlab-editor,
basicstyle=\scriptsize\mlttfamily,
escapechar=`
]
l = -log(eps);
ltol = -log(1e-15);
muvec = sort([linspace(-6,40,60),2]);
c = zeros(size(muvec));
Nval = zeros(size(muvec));
for i=1:length(muvec)
if muvec(i) < 2
 N = @(c,mu) floor(sqrt(2*l*ltol)/pi);
 a = 0;
elseif muvec(i) == 2
 N = @(c,mu) (sqrt(l*ltol)/pi)*...
 sqrt( 1 + (1./c).*(1 + ...
 log(-log( (l-ltol)*(1-c).^2))/ltol));
 a = 1-1/sqrt(l - ltol) + eps;
else
 N = @(c,mu) (sqrt(l*ltol)/pi)*...
 sqrt( 1 + (1./c).*(1 + ...
 ((2 - real(mu))./ltol).*log(1-c)));
 a = 0;
end
c(i) = abs(fminbnd(@(xi) N(xi,muvec(i)),a,1));
Nval(i) = ceil(N(c(i),muvec(i)));
c(i) = abs(fminbnd(@(xi) N(xi,muvec(i)),0,1));
Nval(i) = ceil(real(N(c(i),muvec(i))));
end
\end{lstlisting}
		\end{minipage}
		&
		\begin{minipage}{.46\columnwidth}
			\centering
			\includegraphics[width=0.85\columnwidth]{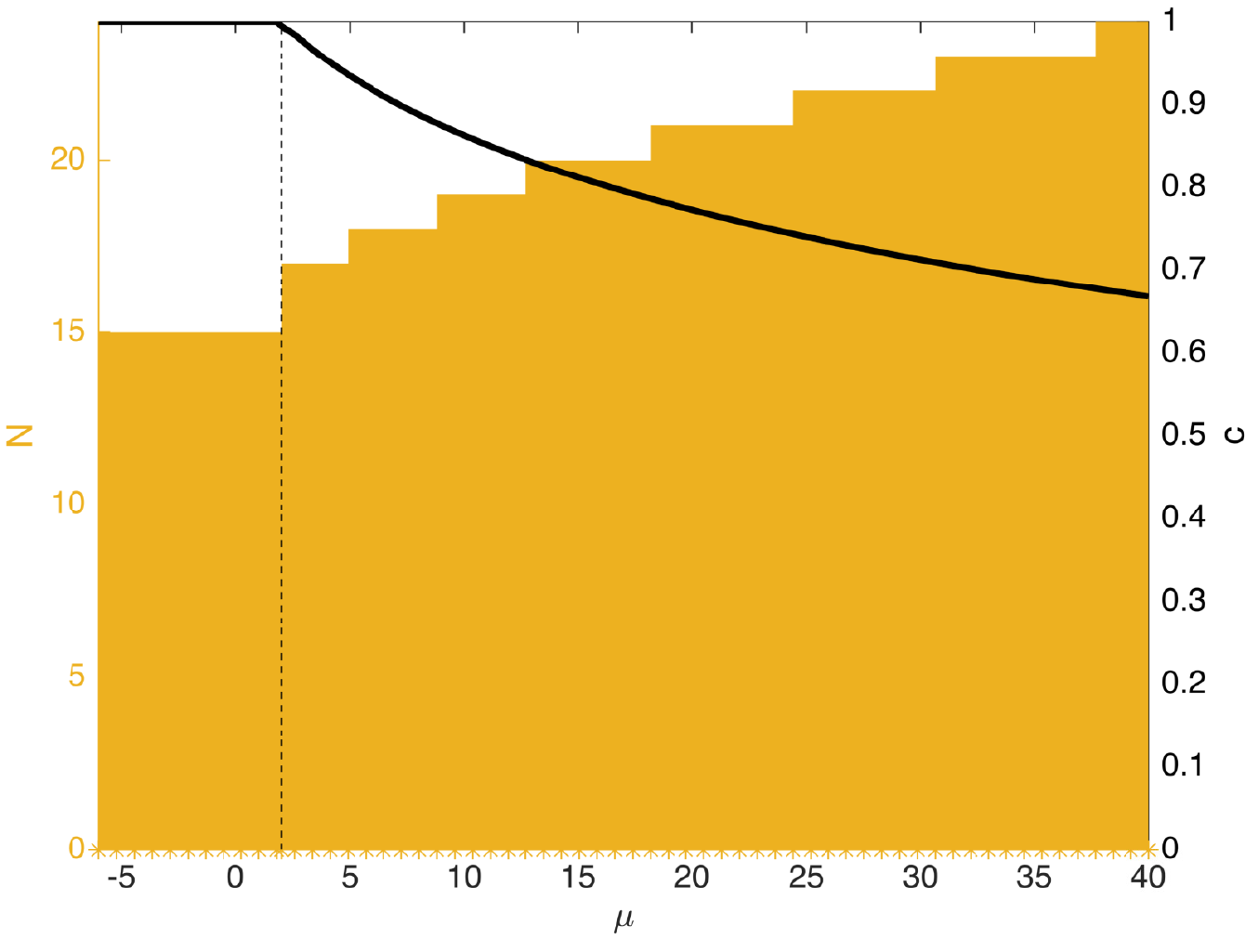}
		\end{minipage}
	\end{tabular}
	\caption{\textsc{MATLAB}\textsuperscript{\textregistered} listing of the optimization for $c$ (left panel), values obtained through the optimization procedure together with the corresponding value of $N$ (right panel) when \emph{double} precision is used, and we pose $\varepsilon = 10^{-15}$ and $\epsilon = 2.2204\times 10^{-16}$. The \emph{dashed} line represents the threshold value $\Re(\mu) = 2$.}
	\label{fig:optimalxi}
\end{figure}
	We observe that as we enlarge the value of $\mu$ the value of $c$ decreases causing an
	increase of the number of quadrature nodes.
\end{remark}

{
	We summarize the previous analysis in the following result.
	\begin{theorem}
		Let $\epsilon$ be the working precision (e.g., single, double, quadruple), and let $\ell= -\log(\epsilon)$. To obtain an absolute error of order $\varepsilon$ on the approximation of the Wright function $t^{\mu-1} W_{\lambda,\mu}(-|x|t^\lambda)$, $\lambda \in (-1,0)$, $\mu \in \mathbb{C}$, $x \in \mathbb{R}$, $t \geq 0$, we select for $\Re(\mu) < 2$
		\[
		N=\left \lfloor{ \frac{\sqrt{2 \ell (-\log \varepsilon)}}{\pi}}\right \rfloor,  \qquad h=\frac{4 \ell}{\pi N^2}, \qquad \gamma=\frac{\pi^2 N^2}{16 t \ell}.
		\]
		We select instead for $\Re(\mu) > 2$
		\[
		N = \left\lfloor \min_{c \in (0,1)}   \frac{\sqrt{\ell (-\log \varepsilon)}}{\pi}\sqrt{ 1 +
			\frac{1}{c}\left( 1+  \frac{2-\Re(\mu)}{(-\log \varepsilon)}  \log(1-c) \right)} \right\rfloor,
		\]
		and use the $\arg\min$ $c$ to compute
		\[
		\xi =  2\left( 1+  \frac{2-\Re(\mu)}{(-\log \varepsilon)}  \log(1-c) \right)^{-1}.
		\]
		While for $\Re(\mu) = 2$
		\[
		N = \left\lfloor \min_{c \in (1 - \nicefrac{1}{\sqrt{\ell+\log(\varepsilon)}},1)} \frac{\sqrt{\ell (-\log \varepsilon)}}{\pi}\sqrt{ 1 +
			\frac{1}{c}\left( 1+  \frac{1}{(-\log \varepsilon)}  \log\left( -\log \left( (\ell+\log \varepsilon) (1-c)^2 \right)\right) \right)} \right\rfloor
		\]
		and use the $\arg\min$ $c$ to compute
		\[
		\xi = 2\left( 1+  \frac{1}{(-\log \varepsilon)}  \log\left(-\log \left( (\ell+\log \varepsilon) (1-c)^2 \right)\right) \right)^{-1}.
		\]
		Finally, for the case $\Re(\mu) \geq 2$ we select
		\[
		h=\frac{(2+\xi c) \ell}{\pi N^2}, \qquad   \gamma=\frac{\pi^2 N^2}{(2+\xi c)^2 t}.
		\]
\end{theorem}}

The method for evaluating the Wright function studied so far is applied to solve some fractional differential problems.
\section{Time-fractional evolution equations on the real line}\label{sec:time-fractional-evolution-equation}

We consider here three time-fractional evolution problems whose solution can be expressed in terms of the Wright function. These model different types of anomalous diffusion, i.e., a diffusion process with a non-linear relationship between the mean squared displacement of the particles being diffused and time. They can be obtained through several generalization of the classical diffusion process leading to the Brownian motion; see, e.g.,~\cite{MR1809268} and {the references}  therein. 

Within this framework, we need to recall the definition of the Caputo fractional derivative
\begin{equation}\label{eq:caputo-derivative}
\frac{\partial^\alpha f(t)}{\partial t^\alpha} = \begin{cases}
f^{(m)}(t), & \text{ if } \alpha = m \in \mathbb{N},\\
\frac{1}{\Gamma(m-\alpha)} \int_{0}^{t} \frac{f^{(m)}(t)}{(t-\tau)^{\alpha+1-m}}\,{\rm d}\tau,& \text{ if } m - 1 < \alpha < m,
\end{cases}
\end{equation}
for $f^{(m)}(t) = \nicefrac{ {\rm d}^m f}{ {\rm d} t^m}$. By means of it, we focus first on the evolution equation~\cite{Mainardi1996,MR1984227,MainardiConsiglio}
\begin{equation}\label{eq:heatproblem}
\frac{\partial^{2\nu} u}{\partial t^{2\nu}} = D \frac{\partial^2 u}{\partial x^2}, \qquad 0 < \nu \leq 1, \quad D > 0,
\end{equation}
for $u(x,t)$ a \emph{causal} function of time, i.e., $u(x,t) \equiv 0$, $\forall\,t<0$. Specifically~\eqref{eq:heatproblem} defines a \emph{fractional diffusion equation} for $\nu \in (0,\nicefrac{1}{2}]$, and a \emph{fractional diffusion-wave equation} for $\nu \in (\nicefrac{1}{2},1]$. We complete such equation to a \emph{Cauchy problem} by posing
\begin{equation}\label{eq:cauchy}
\begin{cases}
u(x,0^+;\nu) = g(x), & x \in \mathbb{R},\\
u(\pm \infty,t;\nu) = 0, & t > 0;
\end{cases}
\end{equation}
whenever $\nu \in (\nicefrac{1}{2},1]$ the initial values of the first time derivative $u_t(x,0^+;\nu) = p(x)$ must be also imposed. The solution of this problem can be expressed by convolutions integrals of the functions in~\eqref{eq:cauchy} with a characteristic function, usually called {Green's} function. This is expressed in terms of the Wright function~\cite{Mainardi1996,Stankovic} (see \eqref{eq:wrightasspecialfunction})
\begin{equation}\label{eq:greenfunctiondefs}
2 \nu |x| \mathcal{G}_C(x,t;\nu) = F_{\nu}(z) = \nu z M_\nu(z) = \nu z W_{-\nu,1-\nu}(-z), \quad z = \nicefrac{|x| t^{-\nu}}{\sqrt{D}},
\end{equation}
where we have denoted the {Green's} function by $\mathcal{G}_C(x,t)$, i.e., the solution for $g(x) = \delta(x)$ the Dirac delta. We then express the general solution of the Cauchy problem \cite[Eq.s (17)-(18)]{MR4309151}~as
\begin{equation}\label{eq:cauchysolution}
u(x,t;\nu) = \begin{cases} \displaystyle \int_{\mathbb{R}} \mathcal{G}_C(\xi,t;\nu) g(x-\xi)\,{\rm d}\xi, & \nu \in(0,\nicefrac{1}{2}],\\[1em]
\displaystyle \int_{\mathbb{R}} \left[ \mathcal{G}_C(\xi,t;\nu)g(x-\xi) + \mathcal{G}_C^{(1)}(\xi,t;\nu)p(x-\xi) \right] \,{\rm d}\xi, &  \nu \in (\nicefrac{1}{2},1],
\end{cases}
\end{equation}
for $\mathcal{G}_C^{(1)}$ a primitive in time of $\mathcal{G}_C$.

Differential equation~\eqref{eq:heatproblem} could be also completed with the following initial and boundary conditions 
\begin{equation}\label{eq:signalling}
\begin{cases}
u(x,0^+;\nu) = 0, & x > 0, \\
u(0^+,t;\nu) = h(t), \; u(+\infty,t,\nu) = 0, & t > 0.
\end{cases}
\end{equation}
This gives us the \emph{{signalling} problem}. To solve it, we introduce  the {Green's} function, i.e., its solution for $h(t) = \delta(t^+)$, 
\begin{equation*}
t \mathcal{G}_S(x,t;\nu) = F_{\nu}(z) = \nu z M_\nu(z) = \nu z W_{-\nu,1-\nu}(-z), \quad z = \nicefrac{|x| t^{-\nu}}{\sqrt{D}}.
\end{equation*}
As for the Cauchy case, to ensure the continuous dependence of the solution on the order parameter $\nu$, we assume $p(x) = 0$ for $\nu \in (\nicefrac{1}{2},1]$, and express it as~\cite[Eq.~(2.6)]{Mainardi1996}
\begin{equation*}
u(x,t;\nu) = \int_{0}^{t} \mathcal{G}_S(x,\tau;\nu) h(t-\tau)\,{\rm d}\tau.
\end{equation*}

The last problem we focus on whose solution can be expressed in terms of the Wright function is the fractional heat conduction in nonhomogeneous media under perfect thermal contact from~\cite{Povstenko1,Povstenko2}. Let us consider two infinite one dimensional {rods joined at}  $x = 0$ with perfect thermal contact, that is, with the same temperature and heat fluxes through the contact point. By using again the Caputo fractional derivative from~\eqref{eq:caputo-derivative} such problem is then stated as
\begin{equation}\label{eq:fheatcond}
\begin{cases}
\displaystyle \frac{\partial^\alpha T_1}{\partial t^\alpha} = a_1 \frac{\partial^2 T_1}{\partial x^2}, & x > 0, \quad 0 < \alpha \leq 2,\\[1em]
\displaystyle \frac{\partial^\beta T_1}{\partial t^\beta} = a_1 \frac{\partial^2 T_1}{\partial x^2}, & x < 0,\quad 0 < \beta \leq 2,\\[1em]
\displaystyle \left.T_1(x,t)\right\rvert_{x = 0^+} = \left.T_2(x,t)\right\rvert_{x = 0^-}, & \text{ continuity of temperatures,}\\[1em]
\displaystyle \left.k_1 D_{\text{RL}}^{1-\alpha} \frac{\partial T_1(x,t)}{\partial x}\right\rvert_{0^+} = \left.k_2 D_{\text{RL}}^{1-\beta} \frac{\partial T_2(x,t)}{\partial x}\right\rvert_{0^-}, & \text{ continuity of heat fluxes,}
\end{cases}
\end{equation}
for $k_1$, $k_2$ the generalized thermal conductivities of the two seminfinite rods and 
\begin{equation*}
D_{\text{RL}}^{\alpha} f(x) = \frac{1}{\Gamma(m-\alpha)} \frac{{\rm d}^m}{{\rm d}t^m} \int_{0}^{t} (t-\tau)^{m-\alpha-1} f(\tau)\,{\rm d}\tau, \quad m-1 < \alpha < m,
\end{equation*}
the Riemann–Liouville derivative of fractional order $\alpha$, which is understood as the Riemann–Liouville fractional integral
\begin{equation*}
I^\alpha f(t) = \frac{1}{\Gamma(\alpha)} \int_{0}^{t} (t-\tau)^{\alpha-1} f(\tau)\,{\rm d}\tau, \quad \alpha > 0,
\end{equation*}
whenever its order turns out to be negative. The solution of~\eqref{eq:fheatcond} with initial condition
\begin{equation*}
t = 0, \quad T_1 = p_0 \delta(x-\rho), \quad x > 0, \quad \rho > 0,
\end{equation*}
for the case $\alpha = \beta$ can be expressed in terms of the Wright function as
\begin{equation}\label{eq:heatrodsolution}
\begin{cases}
\displaystyle T_1(x,t) = \frac{p_0}{2\sqrt{a_1} t^{\nicefrac{\alpha}{2}}} \left[ M_{\frac{\alpha}{2}}\left(\frac{|x-\rho|}{\sqrt{a_1} t^{\nicefrac{\alpha}{2}}}\right) + \frac{\eta-1}{\eta+1} M_{\frac{\alpha}{2}}\left( \frac{x+\rho}{\sqrt{a_1} t^{\nicefrac{\alpha}{2}} } \right)\right], & x \geq 0,\\[1em]
\displaystyle T_2(x,t) = \frac{\eta p_0}{(\eta+1) \sqrt{a_1}t^{\nicefrac{\alpha}{2}}} M_{\frac{\alpha}{2}}\left( \frac{|x|}{\sqrt{a_2}t^{\nicefrac{\alpha}{2}}} + \frac{\rho}{\sqrt{a_1}t^{\nicefrac{\alpha}{2}}}\right), & x \leq 0,
\end{cases}    
\end{equation}
where $\eta = \nicefrac{k_1 \sqrt{a_2}}{k_2 \sqrt{a_1}}$.

\section{Codes and numerical examples} \label{sec:numerexp}

We use this section to numerically validate the results discussed in the sections~\ref{sec:inversion_of_the_Laplace_transform}, \ref{sec:time-fractional-evolution-equation}, and to briefly discuss the implementation of the algorithms for computing the Wright function contained in the repository \href{https://github.com/Cirdans-Home/mwright}{github.com/Cirdans-Home/mwright}. For ease of use, the repository contains a MATLAB\textsuperscript{\textregistered} implementation of the core algorithm in double precision. For better performances and portability, implementations are available in a Fortran module for single, double and quadruple precision.
To test our proposal also in the cases for which alternatives representations of the Wright function are not available, we have included for benchmark purposes a C implementation with arbitrary precision ball arithmetic~\cite{Johansson2017arb} of the series representation~\eqref{wrightfun}.

We make use of these implementation in section~\ref{sec:wrightnumerical} to illustrate the error analysis for the inversion of the Laplace transform. Then, in section~\ref{sec:solving-time-fractional-evolution} we analyze its usage for expressing the solutions of the problems discussed in section~\ref{sec:time-fractional-evolution-equation}.

All the experiments were performed on a Linux machine with an Intel\textsuperscript{\textregistered} Core\texttrademark\, i7-8750H CPU @ 2.20GHz and 16 Gb of memory. The \texttt{m} code are run in MATLAB\textsuperscript{\textregistered} version 9.6.0.1072779 (R2019a). Fortran and C codes are compiled with the \texttt{gnu/9.4.0} suite. To use the single, double and quadruple precision in the Fortran environment we employ the \texttt{iso\_fortran\_env} intrinsic module that provides the needed constants, derived types, and intrinsic procedures.

\subsection{Computing the Wright function}\label{sec:wrightnumerical}

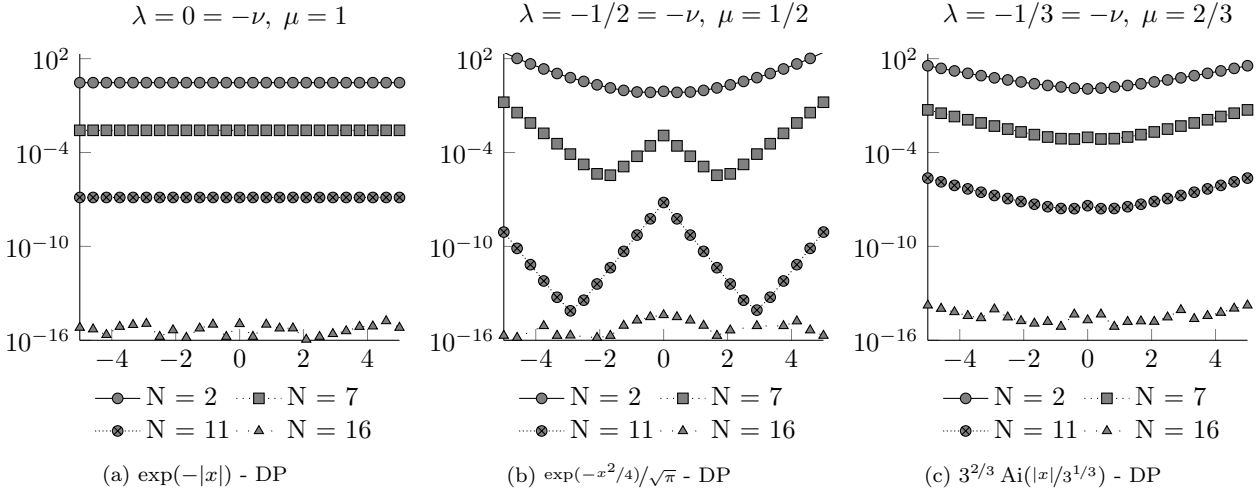
\begin{figure}[t!]
	\centering
	\subfloat[$\exp(-|x|)$ - DP]{% This file was created by matlab2tikz.
\definecolor{mycolor1}{rgb}{0.00000,0.44700,0.74100}%
\definecolor{mycolor2}{rgb}{0.85000,0.32500,0.09800}%
\definecolor{mycolor3}{rgb}{0.92900,0.69400,0.12500}%
\definecolor{mycolor4}{rgb}{0.49400,0.18400,0.55600}%
\pgfplotscreateplotcyclelist{my black white}{%
solid, every mark/.append style={solid, fill=gray}, mark=*\\%
dotted, every mark/.append style={solid, fill=gray}, mark=square*\\%
densely dotted, every mark/.append style={solid, fill=gray}, mark=otimes*\\%
loosely dotted, every mark/.append style={solid, fill=gray}, mark=triangle*\\%
dashed, every mark/.append style={solid, fill=gray},mark=diamond*\\%
loosely dashed, every mark/.append style={solid, fill=gray},mark=*\\%
densely dashed, every mark/.append style={solid, fill=gray},mark=square*\\%
dashdotted, every mark/.append style={solid, fill=gray},mark=otimes*\\%
dashdotdotted, every mark/.append style={solid},mark=star\\%
densely dashdotted,every mark/.append style={solid, fill=gray},mark=diamond*\\%
}
\begin{tikzpicture}

\begin{axis}[%
width=0.25\columnwidth,
height=1.5in,
at={(0.772in,0.481in)},
scale only axis,
xmin=-5,
xmax=5,
ymode=log,
ymin=1e-16,
ymax=200,
yminorticks=true,
axis background/.style={fill=white},
title style={font=\bfseries},
title={$\lambda = 0 = - \nu, \; \mu = 1$},
axis x line*=bottom,
axis y line*=left,
legend columns=2,
legend style={at={(1.6in,-1.5in)},legend cell align=left, align=left, draw=none, fill=none},
cycle list name=my black white
]
\addplot %[color=mycolor1, line width=2.0pt]
  table[row sep=crcr]{%
-5	2.91082385082499\\
-4.58333333333333	2.91082385082499\\
-4.16666666666667	2.91082385082499\\
-3.75	2.91082385082499\\
-3.33333333333333	2.91082385082499\\
-2.91666666666667	2.91082385082499\\
-2.5	2.91082385082499\\
-2.08333333333333	2.91082385082499\\
-1.66666666666667	2.91082385082499\\
-1.25	2.91082385082499\\
-0.833333333333333	2.91082385082499\\
-0.416666666666667	2.91082385082499\\
0	2.91082385082499\\
0.416666666666667	2.91082385082499\\
0.833333333333333	2.91082385082499\\
1.25	2.91082385082499\\
1.66666666666667	2.91082385082499\\
2.08333333333333	2.91082385082499\\
2.5	2.91082385082499\\
2.91666666666667	2.91082385082499\\
3.33333333333333	2.91082385082499\\
3.75	2.91082385082499\\
4.16666666666667	2.91082385082499\\
4.58333333333333	2.91082385082499\\
5	2.91082385082499\\
};
\addlegendentry{N = 2}

\addplot %[color=mycolor2, line width=2.0pt]
  table[row sep=crcr]{%
-5	0.00262514598156318\\
-4.58333333333333	0.0026251459815632\\
-4.16666666666667	0.00262514598156341\\
-3.75	0.00262514598156302\\
-3.33333333333333	0.00262514598156306\\
-2.91666666666667	0.00262514598156339\\
-2.5	0.00262514598156307\\
-2.08333333333333	0.00262514598156333\\
-1.66666666666667	0.00262514598156308\\
-1.25	0.00262514598156318\\
-0.833333333333333	0.00262514598156329\\
-0.416666666666667	0.00262514598156343\\
0	0.00262514598156338\\
0.416666666666667	0.00262514598156343\\
0.833333333333333	0.00262514598156329\\
1.25	0.00262514598156318\\
1.66666666666667	0.00262514598156338\\
2.08333333333333	0.00262514598156344\\
2.5	0.00262514598156307\\
2.91666666666667	0.002625145981563\\
3.33333333333333	0.00262514598156325\\
3.75	0.00262514598156302\\
4.16666666666667	0.0026251459815633\\
4.58333333333333	0.00262514598156337\\
5	0.00262514598156318\\
};
\addlegendentry{N = 7}

\addplot %[color=mycolor3, line width=2.0pt]
  table[row sep=crcr]{%
-5	1.33921290236055e-07\\
-4.58333333333333	1.33921290043536e-07\\
-4.16666666666667	1.33921290358704e-07\\
-3.75	1.33921290235757e-07\\
-3.33333333333333	1.33921290021598e-07\\
-2.91666666666667	1.33921290138038e-07\\
-2.5	1.33921290109844e-07\\
-2.08333333333333	1.33921290377973e-07\\
-1.66666666666667	1.33921290013857e-07\\
-1.25	1.3392129048526e-07\\
-0.833333333333333	1.33921290225154e-07\\
-0.416666666666667	1.3392128991311e-07\\
0	1.33921290057515e-07\\
0.416666666666667	1.3392128991311e-07\\
0.833333333333333	1.33921290225154e-07\\
1.25	1.3392129048526e-07\\
1.66666666666667	1.33921290160809e-07\\
2.08333333333333	1.33921290266518e-07\\
2.5	1.33921290109844e-07\\
2.91666666666667	1.33921290394494e-07\\
3.33333333333333	1.33921290410615e-07\\
3.75	1.33921290235757e-07\\
4.16666666666667	1.33921290134924e-07\\
4.58333333333333	1.33921289873811e-07\\
5	1.33921290236055e-07\\
};
\addlegendentry{N = 11}

\addplot %[color=mycolor4, line width=2.0pt]
  table[row sep=crcr]{%
-5	6.43639478097801e-16\\
-4.58333333333333	5.09175954195695e-16\\
-4.16666666666667	2.23779651284588e-16\\
-3.75	7.37623191694437e-16\\
-3.33333333333333	9.72542355486207e-16\\
-2.91666666666667	1.15405098358923e-15\\
-2.5	1.69066066156624e-16\\
-2.08333333333333	4.45820879994411e-16\\
-1.66666666666667	1.46951618942016e-16\\
-1.25	5.8125886728186e-16\\
-0.833333333333333	1.02183856527064e-15\\
-0.416666666666667	1.68409374933033e-16\\
0	1.11022302462516e-15\\
0.416666666666667	1.68409374933033e-16\\
0.833333333333333	1.02183856527064e-15\\
1.25	5.8125886728186e-16\\
1.66666666666667	5.87806475768064e-16\\
2.08333333333333	1.11455219998603e-16\\
2.5	1.69066066156624e-16\\
2.91666666666667	2.56455774130941e-16\\
3.33333333333333	3.89016942194483e-16\\
3.75	7.37623191694437e-16\\
4.16666666666667	7.83228779496059e-16\\
4.58333333333333	1.69725318065232e-15\\
5	6.43639478097801e-16\\
};
\addlegendentry{N = 16}

\end{axis}
\end{tikzpicture}%
}\hfil
	\subfloat[$\nicefrac{\exp(-\nicefrac{x^2}{4})}{\sqrt{\pi}}$ - DP]{% This file was created by matlab2tikz.
\definecolor{mycolor1}{rgb}{0.00000,0.44700,0.74100}%
\definecolor{mycolor2}{rgb}{0.85000,0.32500,0.09800}%
\definecolor{mycolor3}{rgb}{0.92900,0.69400,0.12500}%
\definecolor{mycolor4}{rgb}{0.49400,0.18400,0.55600}%
\pgfplotscreateplotcyclelist{my black white}{%
solid, every mark/.append style={solid, fill=gray}, mark=*\\%
dotted, every mark/.append style={solid, fill=gray}, mark=square*\\%
densely dotted, every mark/.append style={solid, fill=gray}, mark=otimes*\\%
loosely dotted, every mark/.append style={solid, fill=gray}, mark=triangle*\\%
dashed, every mark/.append style={solid, fill=gray},mark=diamond*\\%
loosely dashed, every mark/.append style={solid, fill=gray},mark=*\\%
densely dashed, every mark/.append style={solid, fill=gray},mark=square*\\%
dashdotted, every mark/.append style={solid, fill=gray},mark=otimes*\\%
dashdotdotted, every mark/.append style={solid},mark=star\\%
densely dashdotted,every mark/.append style={solid, fill=gray},mark=diamond*\\%
}
\begin{tikzpicture}

\begin{axis}[%
width=0.25\columnwidth,
height=1.5in,
at={(0.772in,0.481in)},
scale only axis,
xmin=-5,
xmax=5,
ymode=log,
ymin=1e-16,
ymax=200,
yminorticks=true,
axis background/.style={fill=white},
title style={font=\bfseries},
title={$\lambda = -1/2 = - \nu, \; \mu = 1/2$},
axis x line*=bottom,
axis y line*=left,
legend columns=2,
legend style={at={(1.6in,-1.5in)},legend cell align=left, align=left, draw=none, fill=none},
cycle list name=my black white
]
\addplot %[color=mycolor1, line width=2.0pt]
  table[row sep=crcr]{%
-5	252.98170154567\\
-4.58333333333333	103.384905741602\\
-4.16666666666667	45.7775132697415\\
-3.75	21.8635808081996\\
-3.33333333333333	11.1924497230512\\
-2.91666666666667	6.09278295301352\\
-2.5	3.49936617182174\\
-2.08333333333333	2.11206485448326\\
-1.66666666666667	1.34773655295057\\
-1.25	0.932362166218632\\
-0.833333333333333	0.73504282806833\\
-0.416666666666667	0.698807745823824\\
0	0.813591594927165\\
0.416666666666667	0.698807745823824\\
0.833333333333333	0.73504282806833\\
1.25	0.932362166218632\\
1.66666666666667	1.34773655295057\\
2.08333333333333	2.11206485448325\\
2.5	3.49936617182174\\
2.91666666666667	6.09278295301352\\
3.33333333333333	11.1924497230512\\
3.75	21.8635808081996\\
4.16666666666667	45.7775132697414\\
4.58333333333333	103.384905741602\\
5	252.98170154567\\
};
\addlegendentry{N = 2}

\addplot %[color=mycolor2, line width=2.0pt]
  table[row sep=crcr]{%
-5	0.163460920252024\\
-4.58333333333333	0.0355285259385077\\
-4.16666666666667	0.00772218947515513\\
-3.75	0.00167843307091487\\
-3.33333333333333	0.000364816246541932\\
-2.91666666666667	7.93201122637704e-05\\
-2.5	1.73629495432754e-05\\
-2.08333333333333	4.33783915205041e-06\\
-1.66666666666667	3.53757407533765e-06\\
-1.25	1.27068570500435e-05\\
-0.833333333333333	5.76864493596454e-05\\
-0.416666666666667	0.000265237288666029\\
0	0.00122027661524924\\
0.416666666666667	0.000265237288666029\\
0.833333333333333	5.76864493596454e-05\\
1.25	1.27068570500435e-05\\
1.66666666666667	3.53757407553469e-06\\
2.08333333333333	4.3378391519048e-06\\
2.5	1.73629495432754e-05\\
2.91666666666667	7.93201122641831e-05\\
3.33333333333333	0.000364816246541735\\
3.75	0.00167843307091487\\
4.16666666666667	0.00772218947515512\\
4.58333333333333	0.0355285259385069\\
5	0.163460920252024\\
};
\addlegendentry{N = 7}

\addplot %[color=mycolor3, line width=2.0pt]
  table[row sep=crcr]{%
-5	8.22640011135102e-10\\
-4.58333333333333	7.47496683766992e-11\\
-4.16666666666667	6.79283081843923e-12\\
-3.75	6.16614455483695e-13\\
-3.33333333333333	5.61768113238856e-14\\
-2.91666666666667	7.6335518661364e-15\\
-2.5	3.63782665658642e-14\\
-2.08333333333333	3.95598362678533e-13\\
-1.66666666666667	4.35432842529544e-12\\
-1.25	4.79219662403521e-11\\
-0.833333333333333	5.27391048273201e-10\\
-0.416666666666667	5.80408281064314e-09\\
0	6.38754803890821e-08\\
0.416666666666667	5.80408281064314e-09\\
0.833333333333333	5.27391048273201e-10\\
1.25	4.79219662403521e-11\\
1.66666666666667	4.35432842529544e-12\\
2.08333333333333	3.95307160350466e-13\\
2.5	3.63782665658642e-14\\
2.91666666666667	8.45880071652953e-15\\
3.33333333333333	5.63746169975613e-14\\
3.75	6.16614455483695e-13\\
4.16666666666667	6.79224101963235e-12\\
4.58333333333333	7.47490814211599e-11\\
5	8.22640011135102e-10\\
};
\addlegendentry{N = 11}

\addplot %[color=mycolor4, line width=2.0pt]
  table[row sep=crcr]{%
-5	1.99092874543732e-16\\
-4.58333333333333	1.46738884851571e-16\\
-4.16666666666667	0\\
-3.75	8.27392761467554e-16\\
-3.33333333333333	1.97805673675653e-16\\
-2.91666666666667	2.06312212598281e-16\\
-2.5	0\\
-2.08333333333333	1.45601164033321e-16\\
-1.66666666666667	1.97037351250981e-16\\
-1.25	8.72474420186105e-16\\
-0.833333333333333	1.87272281180645e-15\\
-0.416666666666667	3.0826627324971e-15\\
0	4.13242005825774e-15\\
0.416666666666667	3.0826627324971e-15\\
0.833333333333333	1.87272281180645e-15\\
1.25	8.72474420186105e-16\\
1.66666666666667	1.97037351250981e-16\\
2.08333333333333	4.36803492099963e-16\\
2.5	0\\
2.91666666666667	8.25248850393125e-16\\
3.33333333333333	0\\
3.75	8.27392761467554e-16\\
4.16666666666667	1.65143665922532e-15\\
4.58333333333333	4.40216654554713e-16\\
5	1.99092874543732e-16\\
};
\addlegendentry{N = 16}

\end{axis}

\end{tikzpicture}%
}\hfil
	\subfloat[$3^{\nicefrac{2}{3}} \operatorname{Ai}(\nicefrac{|x|}{3^{\nicefrac{1}{3}}})$ - DP]{% This file was created by matlab2tikz.
\definecolor{mycolor1}{rgb}{0.00000,0.44700,0.74100}%
\definecolor{mycolor2}{rgb}{0.85000,0.32500,0.09800}%
\definecolor{mycolor3}{rgb}{0.92900,0.69400,0.12500}%
\definecolor{mycolor4}{rgb}{0.49400,0.18400,0.55600}%
\pgfplotscreateplotcyclelist{my black white}{%
solid, every mark/.append style={solid, fill=gray}, mark=*\\%
dotted, every mark/.append style={solid, fill=gray}, mark=square*\\%
densely dotted, every mark/.append style={solid, fill=gray}, mark=otimes*\\%
loosely dotted, every mark/.append style={solid, fill=gray}, mark=triangle*\\%
dashed, every mark/.append style={solid, fill=gray},mark=diamond*\\%
loosely dashed, every mark/.append style={solid, fill=gray},mark=*\\%
densely dashed, every mark/.append style={solid, fill=gray},mark=square*\\%
dashdotted, every mark/.append style={solid, fill=gray},mark=otimes*\\%
dashdotdotted, every mark/.append style={solid},mark=star\\%
densely dashdotted,every mark/.append style={solid, fill=gray},mark=diamond*\\%
}
\begin{tikzpicture}

\begin{axis}[%
width=0.25\columnwidth,
height=1.5in,
at={(0.772in,0.481in)},
scale only axis,
xmin=-5,
xmax=5,
ymode=log,
ymin=1e-16,
ymax=200,
yminorticks=true,
axis background/.style={fill=white},
title style={font=\bfseries},
title={$\lambda = -1/3 = - \nu, \; \mu = 2/3$},
axis x line*=bottom,
axis y line*=left,
legend columns=2,
legend style={at={(1.6in,-1.5in)},legend cell align=left, align=left, draw=none, fill=none},
cycle list name=my black white
]
\addplot %[color=mycolor1, line width=2.0pt]
  table[row sep=crcr]{%
-5	35.100946632256\\
-4.58333333333333	23.7866581598429\\
-4.16666666666667	16.392740189559\\
-3.75	11.4835249861634\\
-3.33333333333333	8.17357036487473\\
-2.91666666666667	5.90936599476757\\
-2.5	4.34035005831199\\
-2.08333333333333	3.24177922112168\\
-1.66666666666667	2.46808903673648\\
-1.25	1.92443655696228\\
-0.833333333333333	1.54928300343777\\
-0.416666666666667	1.30388590671308\\
0	1.16641429145616\\
0.416666666666667	1.30388590671308\\
0.833333333333333	1.54928300343777\\
1.25	1.92443655696228\\
1.66666666666667	2.46808903673648\\
2.08333333333333	3.24177922112168\\
2.5	4.34035005831199\\
2.91666666666667	5.90936599476757\\
3.33333333333333	8.17357036487474\\
3.75	11.4835249861634\\
4.16666666666667	16.3927401895589\\
4.58333333333333	23.7866581598429\\
5	35.100946632256\\
};
\addlegendentry{N = 2}

\addplot %[color=mycolor2, line width=2.0pt]
  table[row sep=crcr]{%
-5	0.0529821544727944\\
-4.58333333333333	0.0324306082365764\\
-4.16666666666667	0.0199706274781328\\
-3.75	0.0123897646799814\\
-3.33333333333333	0.00775718364204686\\
-2.91666666666667	0.00491153498681992\\
-2.5	0.00315373986713108\\
-2.08333333333333	0.00206307440133608\\
-1.66666666666667	0.00138769299387815\\
-1.25	0.000980804704677108\\
-0.833333333333333	0.00076661109081967\\
-0.416666666666667	0.000729934291244036\\
0	0.000934911401088799\\
0.416666666666667	0.000729934291244036\\
0.833333333333333	0.00076661109081967\\
1.25	0.000980804704677108\\
1.66666666666667	0.00138769299387815\\
2.08333333333333	0.00206307440133504\\
2.5	0.00315373986713108\\
2.91666666666667	0.00491153498681933\\
3.33333333333333	0.00775718364204719\\
3.75	0.0123897646799814\\
4.16666666666667	0.0199706274781321\\
4.58333333333333	0.032430608236577\\
5	0.0529821544727944\\
};
\addlegendentry{N = 7}

\addplot %[color=mycolor3, line width=2.0pt]
  table[row sep=crcr]{%
-5	2.31045314878271e-06\\
-4.58333333333333	1.34261713955476e-06\\
-4.16666666666667	7.91734742221577e-07\\
-3.75	4.74372299556711e-07\\
-3.33333333333333	2.89173940749841e-07\\
-2.91666666666667	1.79620173570332e-07\\
-2.5	1.13890793951548e-07\\
-2.08333333333333	7.3905605679551e-08\\
-1.66666666666667	4.93403695623865e-08\\
-1.25	3.44151142215259e-08\\
-0.833333333333333	2.64144924446848e-08\\
-0.416666666666667	2.5857498479363e-08\\
0	3.99165205967557e-08\\
0.416666666666667	2.5857498479363e-08\\
0.833333333333333	2.64144924446848e-08\\
1.25	3.44151142215259e-08\\
1.66666666666667	4.93403696813682e-08\\
2.08333333333333	7.39056046441255e-08\\
2.5	1.13890793951548e-07\\
2.91666666666667	1.79620172779777e-07\\
3.33333333333333	2.89173940749841e-07\\
3.75	4.74372299556711e-07\\
4.16666666666667	7.91734741599373e-07\\
4.58333333333333	1.34261713937952e-06\\
5	2.31045314878271e-06\\
};
\addlegendentry{N = 11}

\addplot %[color=mycolor4, line width=2.0pt]
  table[row sep=crcr]{%
-5	1.71013674979627e-14\\
-4.58333333333333	1.05149515223987e-14\\
-4.16666666666667	6.42942539629809e-15\\
-3.75	3.6404337575709e-15\\
-3.33333333333333	2.6449193203169e-15\\
-2.91666666666667	9.68430160473201e-15\\
-2.5	2.96285708040835e-15\\
-2.08333333333333	1.72570892663818e-15\\
-1.66666666666667	1.18981680996062e-15\\
-1.25	1.52288783392817e-15\\
-0.833333333333333	7.46298548559209e-16\\
-0.416666666666667	4.55000603049022e-15\\
0	1.95438478873191e-15\\
0.416666666666667	4.55000603049022e-15\\
0.833333333333333	7.46298548559209e-16\\
1.25	1.52288783392817e-15\\
1.66666666666667	1.66574353394488e-15\\
2.08333333333333	1.55313803397436e-15\\
2.5	2.96285708040835e-15\\
2.91666666666667	8.69610756343283e-15\\
3.33333333333333	2.33375234145609e-15\\
3.75	3.6404337575709e-15\\
4.16666666666667	5.80722293859181e-15\\
4.58333333333333	1.08654499064787e-14\\
5	1.71013674979627e-14\\
};
\addlegendentry{N = 16}

\end{axis}
\end{tikzpicture}%
}
	
	\caption{Relative error with respect to the number of quadrature points~$N$ for cases in which we have an exact alternative representations of the Mainardi function $M_\nu(|x|) = W_{-\nu,1-\nu}(-|x|)$ in~\eqref{eq:special-cases}. Results in double precision (DP) are computed with the Matlab routine.}
	\label{fig:error}
\end{figure}

We consider first the Mainardi functions $M_\nu(|x|)$, that are the cases of the form $W_{\lambda,\mu}(-|x|)$, $\lambda = -\nu$, $\mu = 1-\nu$, for which we have a closed form relation~\eqref{eq:special-cases}. For these cases we report in Figure~\ref{fig:error} the convergence with respect to the number of nodes $N$ of the approximation in~\eqref{eq:trapN} to these close form expressions by looking at the relative error with respect to the exact representation. The last value of $N$ is the one given in~\eqref{eq:parmin2} from which we obtain the expected value. 
For the two cases in which the function can be represented as an exponential we can also test the routines for calculating the quadruple precision, see Figure~\ref{fig:error-quadruple}.
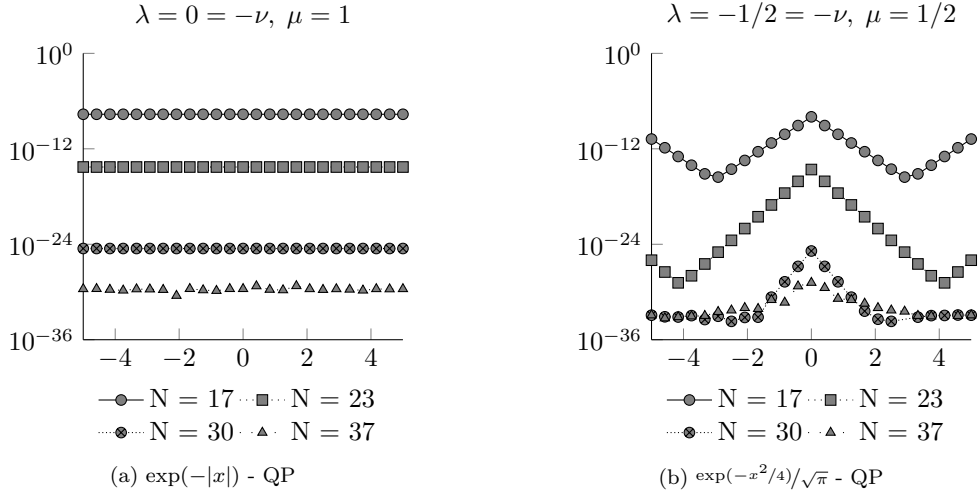
\begin{figure}[htbp]
	\centering
	\subfloat[$\exp(-|x|)$ - QP]{\pgfplotscreateplotcyclelist{my black white}{%
solid, every mark/.append style={solid, fill=gray}, mark=*\\%
dotted, every mark/.append style={solid, fill=gray}, mark=square*\\%
densely dotted, every mark/.append style={solid, fill=gray}, mark=otimes*\\%
loosely dotted, every mark/.append style={solid, fill=gray}, mark=triangle*\\%
dashed, every mark/.append style={solid, fill=gray},mark=diamond*\\%
loosely dashed, every mark/.append style={solid, fill=gray},mark=*\\%
densely dashed, every mark/.append style={solid, fill=gray},mark=square*\\%
dashdotted, every mark/.append style={solid, fill=gray},mark=otimes*\\%
dashdotdotted, every mark/.append style={solid},mark=star\\%
densely dashdotted,every mark/.append style={solid, fill=gray},mark=diamond*\\%
}
\begin{tikzpicture}
\begin{axis}[%
width=0.25\columnwidth,
height=1.5in,
at={(0.772in,0.481in)},
scale only axis,
xmin=-5,
xmax=5,
ymode=log,
ymin=1e-36,
ymax=1,
yminorticks=true,
axis background/.style={fill=white},
title style={font=\bfseries},
title={$\lambda = 0 = - \nu, \; \mu = 1$},
axis x line*=bottom,
axis y line*=left,
legend columns=2,
legend style={at={(1.6in,-1.7in)},legend cell align=left, align=left, draw=none, fill=none},
cycle list name=my black white,
%every axis plot/.append style={ultra thick}
]

\addplot table[col sep=comma] {mainardi-m0-quadruple-17.out};
\addlegendentry{N = 17}

\addplot table[col sep=comma] {mainardi-m0-quadruple-23.out};
\addlegendentry{N = 23}

\addplot table[col sep=comma] {mainardi-m0-quadruple-30.out};
\addlegendentry{N = 30}

\addplot table[col sep=comma] {mainardi-m0-quadruple-37.out};
\addlegendentry{N = 37}

\end{axis}
\end{tikzpicture}}\hfil
	\subfloat[$\nicefrac{\exp(-\nicefrac{x^2}{4})}{\sqrt{\pi}}$ - QP]{\pgfplotscreateplotcyclelist{my black white}{%
solid, every mark/.append style={solid, fill=gray}, mark=*\\%
dotted, every mark/.append style={solid, fill=gray}, mark=square*\\%
densely dotted, every mark/.append style={solid, fill=gray}, mark=otimes*\\%
loosely dotted, every mark/.append style={solid, fill=gray}, mark=triangle*\\%
dashed, every mark/.append style={solid, fill=gray},mark=diamond*\\%
loosely dashed, every mark/.append style={solid, fill=gray},mark=*\\%
densely dashed, every mark/.append style={solid, fill=gray},mark=square*\\%
dashdotted, every mark/.append style={solid, fill=gray},mark=otimes*\\%
dashdotdotted, every mark/.append style={solid},mark=star\\%
densely dashdotted,every mark/.append style={solid, fill=gray},mark=diamond*\\%
}
\begin{tikzpicture}
\begin{axis}[%
width=0.25\columnwidth,
height=1.5in,
at={(0.772in,0.481in)},
scale only axis,
xmin=-5,
xmax=5,
ymode=log,
ymin=1e-36,
ymax=1,
yminorticks=true,
axis background/.style={fill=white},
title style={font=\bfseries},
title={$\lambda = -1/2 = - \nu, \; \mu = 1/2$},
axis x line*=bottom,
axis y line*=left,
legend columns=2,
legend style={at={(1.6in,-1.7in)},legend cell align=left, align=left, draw=none, fill=none},
cycle list name=my black white,
%every axis plot/.append style={ultra thick}
]

\addplot table[col sep=comma] {mainardi-m12-quadruple-17.out};
\addlegendentry{N = 17}

\addplot table[col sep=comma] {mainardi-m12-quadruple-23.out};
\addlegendentry{N = 23}

\addplot table[col sep=comma] {mainardi-m12-quadruple-30.out};
\addlegendentry{N = 30}

\addplot table[col sep=comma] {mainardi-m12-quadruple-37.out};
\addlegendentry{N = 37}

\end{axis}
\end{tikzpicture}}
	
	\caption{Relative error with respect to the number of quadrature points~$N$ for cases in which we have an exact alternative representations of the Mainardi function $M_\nu(|x|) = W_{-\nu,1-\nu}(-|x|)$ in~\eqref{eq:special-cases}. Results in quadruple precision (QP) are computed with the Fortran routine.}
	\label{fig:error-quadruple}
\end{figure}
To evaluate the accuracy in the general case, we make use of the extended precision ball arithmetic library~\cite{Johansson2017arb}. For this case we perform a comparison also for the case in which the parameter $\mu \in \mathbb{C}$. 
We consider only the number of quadrature points determined by~\eqref{eq:valueofN}, and report the results in the \emph{heatmap}s in Figure~\ref{fig:complexmu}. 
\begin{figure}[htbp]
	\centering
	\includegraphics[width=\columnwidth]{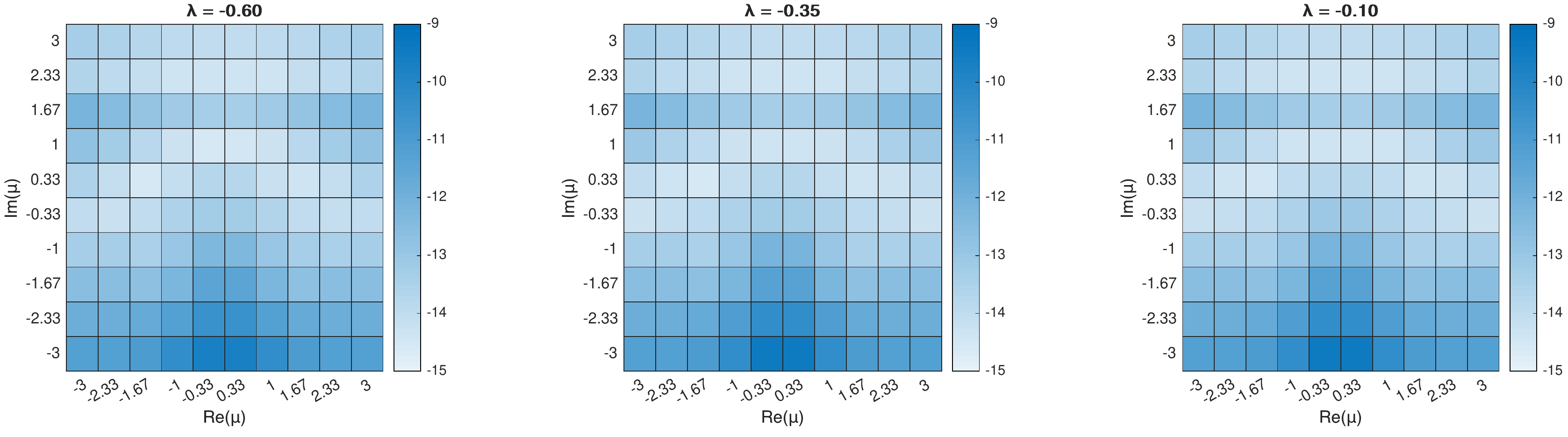}
	\caption{Heatmaps for the relative error (color map is in log-scale) between the values computed with the series computed to the $1000$th entry in augmented 128 bits precision and the values computed with our algorithm for the cases in which a complex $\mu$ is used. The number of quadrature points has been fixed as in the analysis in section~\ref{sec:selecting_the_quadrature_nodes}.}
	\label{fig:complexmu}
\end{figure}
Each of the box in the plots consider the norm-wise relative error in $2$-norm over the interval $[-5,0]$, the color map is in logarithmic scale. In all the cases we reach at least a relative error of the order of $10^{-10}$ for the \emph{double} precision routine.

\message{The column width is: \the\columnwidth}

\subsection{Solving the time-fractional evolution equation}\label{sec:solving-time-fractional-evolution}

In this section, we deal with the approximation of the exact solution, represented by convolution with {the} Green's function, of the differential problems involving fractional derivatives in time discussed in section~\ref{sec:time-fractional-evolution-equation}. 

We focus first on the analysis of the solution of the Cauchy problem in~\eqref{eq:cauchysolution}. To calculate the solution we need to perform a convolution between the values obtained with the computation of the Wright function and the initial condition through a fast Fourier transform.
\begin{figure}[htbp]
	\centering
	\subfloat{\includegraphics{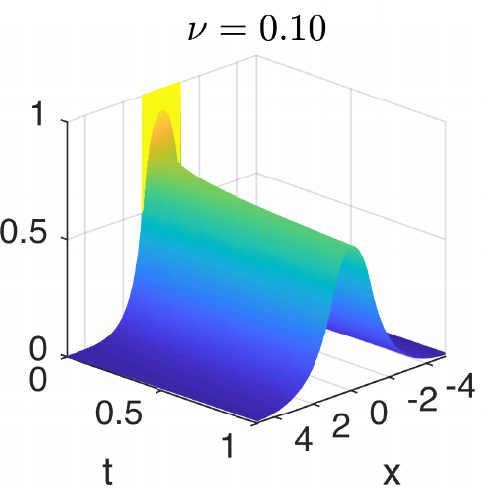}}\hfil
	\subfloat{\includegraphics{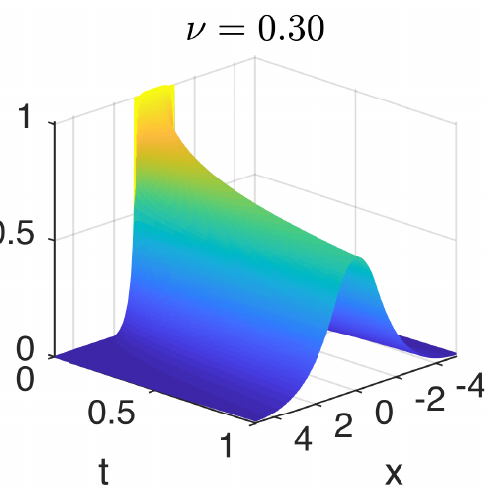}}\hfil
	\subfloat{\includegraphics{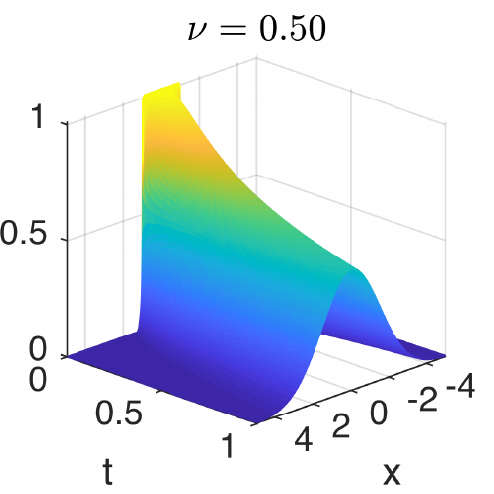}}
	\caption{Cauchy problem. Solution of the Cauchy problem as in~\eqref{eq:cauchysolution} were the initial condition $g(x)$ is approximated by a square pulse centered at zero and amplitude $1$, for several values of the parameter $\nu \in (0,\nicefrac{1}{2}]$, and $D = 1$.}
	\label{fig:cauchy_problem}
\end{figure}
In Figure~\ref{fig:cauchy_problem} we use the idea contained in the following snippet of MATLAB\textsuperscript{\textregistered} code
\begin{lstlisting}[style=Matlab-editor]
xlr=5; np=256; dx=2*xlr/np; nu=0.5;
x = linspace(-xlr,xlr-dx,np);
g = real((abs(x) <= 1));
t = 1;
Gc = mwright(x,t,-nu,1-nu)/2;
fg = fft(g);
fGc = fft(Gc);
u = dx*fftshift(ifft( fg.*fGc ));
plot(x,g,'k',x,u,'r');
\end{lstlisting}
to produce the solution for the case of $\nu \in (0,\nicefrac{1}{2}]$. The expected sub-diffusion behavior is qualitatively observed with respect to the case of the classical solution ($\nu = \nicefrac{1}{2}$). That is, the convergence of the solution to the zero stationary state is slower, the more the order $\nu$ approaches 0. 

We observe that to obtain an approximate solution using a finite difference method, rather than a method based on a weak formulation, it would be necessary to set a domain truncation and study the decay to zero of the function to keep the error under control while also having to solve the associated linear system. Moreover, methods for marching in time fractional differential equations have usually a low approximation order and needs to retain weighty memory terms, see, e.g., the general discussion in~\cite{MR2176361}. To illustrate this behavior we consider here the usage of the fractional trapezoidal rule~\cite{MR3327641} that is an example of a Fractional Linear Multistep Method (FLMM). We discretize the spatial term of~\eqref{eq:heatproblem} by second order centered finite differences over the domain $[-20,20]$ so that is reasonable to close the system of equations by using homogeneous Dirichlet boundary conditions. 
\begin{figure}[htbp]
	\centering
	\subfloat[Comparison of the obtained solution, absolute error ($h = 0.01$, $n = 256$)\label{fig:flmm-1}]{\input{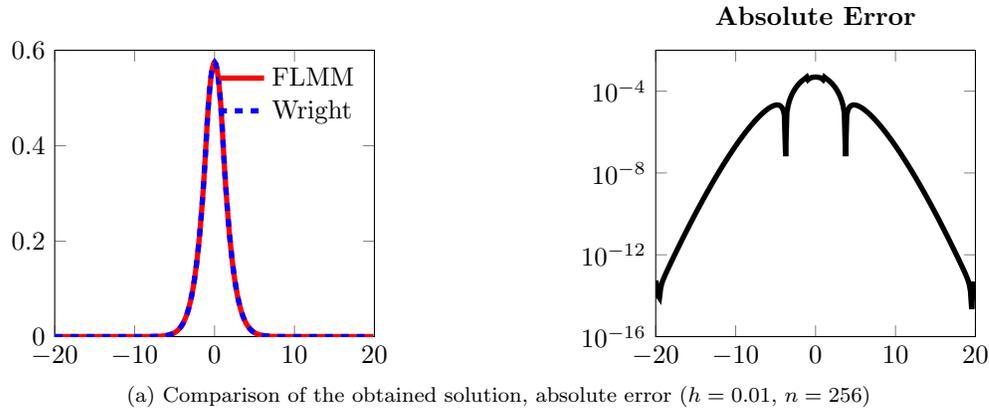}
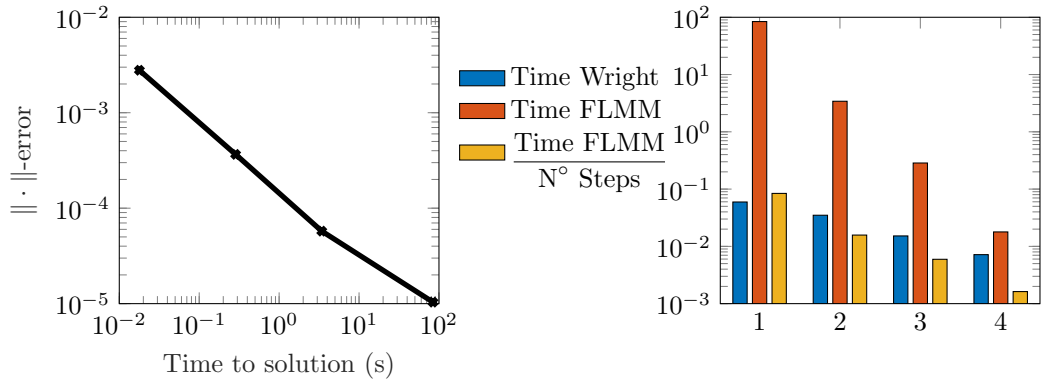}
	
	\subfloat[Time versus error to reach the last time step ($t = 1$) and bar plots comparing the time needed by the Wright function to compute the solution on the last time step and compare it with both the time needed to march the solution to $t=1$ and the \emph{average} time-per-step\label{fig:flmm-2}]{% This file was created by matlab2tikz.
\definecolor{mycolor1}{rgb}{0.00000,0.44700,0.74100}%
\definecolor{mycolor2}{rgb}{0.85000,0.32500,0.09800}%
\definecolor{mycolor3}{rgb}{0.92900,0.69400,0.12500}%
\begin{tikzpicture}

\begin{axis}[%
width=0.25\columnwidth,
height=1.5in,
at={(0in,0in)},
scale only axis,
xmode=log,
xmin=0.01,
xmax=100,
xminorticks=true,
xlabel style={font=\color{white!15!black}},
xlabel={Time to solution (s)},
ymode=log,
ymin=1e-05,
ymax=0.01,
yminorticks=true,
ylabel style={font=\color{white!15!black}},
ylabel={$\|\cdot\|$-error},
axis background/.style={fill=white},
legend style={legend cell align=left, align=left, draw=white!15!black}
]
\addplot [color=black, line width=2.0pt, mark=x, mark options={solid, black}]
  table[row sep=crcr]{%
84.170734	1.03704358421775e-05\\
3.417063	5.74907072805608e-05\\
0.285506	0.000366035365327366\\
0.017868	0.00278448453197216\\
};
%\addlegendentry{data1}

\end{axis}

\begin{axis}[%
width=0.25\columnwidth,
height=1.5in,
at={(0.47\columnwidth,0in)},
scale only axis,
bar shift auto,
log origin=infty,
xmin=0.511111111111111,
xmax=4.48888888888889,
xtick={1, 2, 3, 4},
ymode=log,
ymin=0.001,
ymax=100,
yminorticks=true,
axis background/.style={fill=white},
legend style={at={(-0.03\columnwidth,0.397in)},legend cell align=left, align=left, draw=none, fill=none}
]
\addplot[ybar, bar width=0.178, fill=mycolor1, draw=black, area legend] table[row sep=crcr] {%
1	0.059535\\
2	0.034955\\
3	0.015214\\
4	0.007172\\
};
\addlegendentry{Time Wright}

\addplot[ybar, bar width=0.178, fill=mycolor2, draw=black, area legend] table[row sep=crcr] {%
1	84.170734\\
2	3.417063\\
3	0.285506\\
4	0.017868\\
};
\addlegendentry{Time FLMM}

\addplot[ybar, bar width=0.178, fill=mycolor3, draw=black, area legend] table[row sep=crcr] {%
1	0.0840866473526473\\
2	0.0157468341013825\\
3	0.00594804166666667\\
4	0.00162436363636364\\
};
\addlegendentry{$\displaystyle \frac{\text{Time FLMM}}{\displaystyle \text{N}^\circ \text{ Steps} }$}

\end{axis}

\end{tikzpicture}%
}
	
	\caption{Comparison of the Trapezoidal FLMM for the Cauhcy problem with $\nu = 0.3$, $D=1$.}
\end{figure}
In Figure~\ref{fig:flmm-1} we report a visual comparison between the solution obtained by using the finite difference method on $n = 256$ grid points (continuous line) and the one obtained by means of the representation~\eqref{eq:cauchysolution} (dashed line) on the same spatial grid. Having selected a time step for the FLMM method of $h =10^{-2}$ we report on the second panel the absolute error across the spatial domain on the last time step ($t = 1$). In Figure~\ref{fig:flmm-2} we have a time versus infinity norm error graph in which the solid curve represents the solution time of the whole procedure on four logarithmically equispaced $h$ values between $10^{-3}$ and $10^{-1}$ and $n = 2048,1024,512,256$. The bar graph depicts instead the time needed by the Wright function to compute the solution on the last time step and compare it with both the time needed to march the solution to $t=1$ and the \emph{average} time-per-step obtained as the total time divided by the number of time-steps. What we observe is that the advantage is more evident when greater accuracy is required, and in general the computation of a single instance of the Wright function is comparable in these cases with the average time required for a single time-step. Observe also that the average value that we report for completeness is not entirely indicative, because it amortizes in the same way between the steps the time due to the treatment of the integration queue, which is instead an increasing quantity with the number of steps.

\begin{figure}[htbp]
	\centering
	\subfloat[$p_0 = 1$, $\rho = 0.5$, $a_1 = 3$, $a_2 = 1$, $k_1 = 2$, $k_2 = 6$]{\includegraphics[width=\columnwidth]{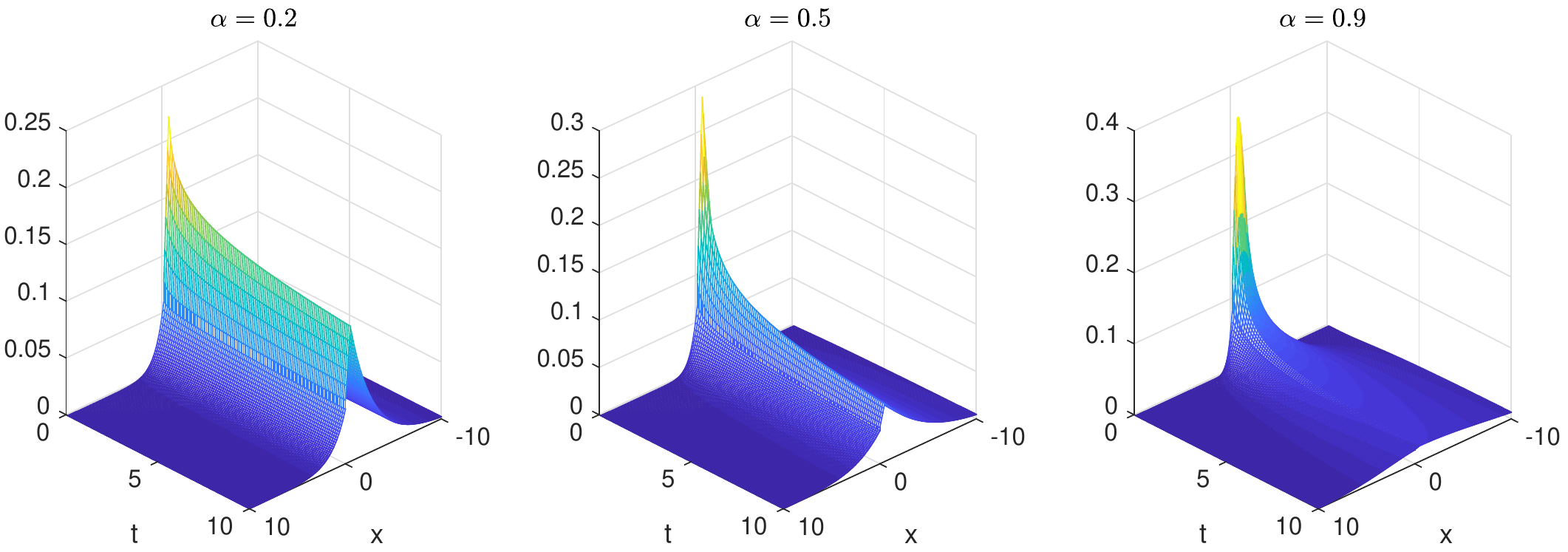}}
	\caption{Solution of the fractional heat conduction in nonhomogeneous media under perfect thermal contact for different parameters choices.}
	\label{fig:heatconduction}
\end{figure}

The other numerical experiment we perform is the heat conduction problem in~\eqref{eq:fheatcond}. In Figure~\ref{fig:heatconduction} we draw the solution~\eqref{eq:heatrodsolution} for the same physical parameters, and different values of the fractional order $\alpha$ computed by means of the MATLAB\textsuperscript{\textregistered} version of the Wright function. 
We observe again that also in this case, wanting to solve the problem with a discretization approach, we would face the same problem of having to discuss the truncation of the domain together with that of having to integrate the two fractional differential equations in time. It should also be noted that in this case the coupling conditions are two fractional integrals, this makes the associated linear systems slightly more complex than those of the previous case.

\section{Conclusions} \label{sec:conclusion}
In this work we have proposed an algorithm  for the computation of the Wright function of the second type based on the inversion of the Laplace transform. As far as we know, this is the first available and open-source code to effectively evaluate this function. The proposed method appears to be also a valid tool for solving fractional diffusion-wave equations thanks to its modest computational cost and high accuracy.

\section*{Acknowledgements} 
This work was partially supported by GNCS-INdAM. The authors are members of the INdAM research group
GNCS. {The authors thank the anonymous referees for the valuable suggestions that led to improve the paper}.

\bibliographystyle{abbrv}
\bibliography{bibliografia}

\end{document}